\documentclass[reqno]{amsart}
\usepackage{amssymb}
\usepackage{graphicx}
 
\def\part{{\partial}}

\def\AA{\mathcal A}

\def\A{\bold A}
\def\L{\mathcal L}
\def\C{\bbb C}
\def\CC{\bold Z}

\def\ZZ{\bold Z}

\def\O{\mathcal O}
\def\C{\Bbb C}

\def\B{\bold B}
\def\BB{\mathcal B}
\def\R{{\Bbb R}}

\newtheorem{Theorem}{Theorem}[section]
\newtheorem{Lemma}{Lemma}[section]
\newtheorem{Proposition}{Proposition}[section]
\newtheorem{Corollary}{Corollary}

\newtheorem{Example}{Example}[section]

\begin{document}

\title[{\sc Flags in Zero Dimensional Complete Intersections and Indices}]
{Flags in Zero Dimensional Complete Intersection Algebras and Indices of Real Vector
Fields \\
{\footnotesize }}
\author{L. Giraldo, X. G\'omez-Mont and P. Marde\v si\'c}
\address{ \hskip -4mm - Dep. Matem\'{a}ticas,  Universidad de C\'adiz \& (current
address) Dep. Geometr\'{\i}a y Topolog\'{\i}a,  F. Matem\'aticas,
Universidad Complutense de Madrid, 28040 Madrid,
 Spain;
\newline - CIMAT, A.P. 402 Guanajuato, 36000
M\'exico;
\newline - Institut de Math\'{e}matique de Bourgogne, U.M.R.
5584 du C.N.R.S.
\newline
{Universit\'{e} de Bourgogne, B.P. 47 870}\\
{21078 Dijon Cedex, France }} \email{luis.giraldo@uca.es;
gmont@cimat.mx; mardesic@u-bourgogne.fr}

\thanks{Partially supported by Plan Nacional I+D grant no.
MTM2004-07203-C02-02, Spain and CONACYT 40329, M\'exico}
\date{\today}
\begin{abstract}
We introduce  bilinear forms in a flag in a complete intersection
local $\R$-algebra of dimension 0, related to the Eisenbud-Levine,
Khimshiashvili bilinear form.
We give a variational interpretation of these forms in terms
of Jantzen's filtration and bilinear forms.
 We use the signatures of these forms to compute in the real case
the constant relating the GSV-index with the signature function of
vector fields tangent to an even dimensional hypersurface
singularity, one being topologically defined and the other
computable with finite dimensional commutative algebra methods.
\end{abstract}
\maketitle


\section{Introduction}
\vskip 1cm Let $f_1,\ldots,f_n:\R^{n} \rightarrow \R$ be  germs of
real analytic functions that form a regular sequence as
holomorphic functions and let
\begin{equation}
\label{d1} \A:=\frac{{\mathcal A}_{\R^{n},0}} {
(f_1,\ldots,f_{n})}
\end{equation} be the
quotient finite dimensional algebra, where
${\mathcal A}_{\R^{n},0}$ is the algebra of germs of real analytic functions
on $\R^{n}$ with coordinates $x_1,\ldots,x_{n}$. The class of the
Jacobian
\begin{equation}\label{d2}
J= det\left(\frac{\partial f_i}{\partial
x_j}\right)_{i,j=1,\ldots,n} \hskip 1cm , \hskip 1cm
J_\A := [J]_\A \in \A
\end{equation}
generates the socle (the unique minimal
non-zero ideal) of the algebra $\A$. A symmetric bilinear form
\begin{equation}\label{d3}
<\  , \,>_{L_\A}: \A \times \A {\stackrel{\cdot}\rightarrow} \A
{\stackrel{L_\A }\rightarrow} \R
\end{equation}
is defined by composing multiplication in $\A$ with any
linear map $L_\A:\A \longrightarrow \R$ sending $J_\A$ to a
positive number. The theory of Eisenbud-Levine and Khimshiashvili
asserts that this bilinear form is nondegenerate and that its
signature $\sigma_\A$ is independent of the choice of $L_\A$
(see \cite{EL}, \cite{K}).

\vskip 3mm Let $f \in \A$ be an element in the maximal ideal. We
define a flag of ideals in $\A$:
\begin{equation}\label{a}
 K_m := Ann_\A(f) \cap (f^{m-1}), \hskip 5mm m\geq1, \hskip 5mm
 0 \subset  K_{\ell +1} \subset \cdots \subset K_1 \subset K_0 := \A
\end{equation}
 and a family of bilinear forms
\begin{equation}\label{aaa}
<\  ,\ >_{L_\A,f,m}:   K_{m} \times  K_{m} \rightarrow  \R
\hskip 3mm,\hskip 3mm <a,a^\prime>_{L_\A,f,m} =
<\frac{a}{f^{m-1}},a^\prime>_{L_\A}, \end{equation} defined for
$m=0,\ldots,\ell+1.$ The division by $f^{m-1}$ is defined up to
elements in $Ann_\A(f^{m-1})$, but as $a'\in (f^{m-1})$, the
last expression in (\ref{aaa}) is well defined. We call the form
$<\ ,\ >_{L_\A,f,m}$, the order $m$ bilinear form on the algebra
$\A$, with respect to $f$. In Section 1 we prove:

 \vskip 3mm \begin{Theorem}\label{th1} For $m=0,\ldots,\ell+1$ the order $m$
bilinear form $<  ,\,>_{L_\A,f,m}$ on $K_{m}$  induces a
non-degenerate bilinear form

\begin{equation}\label{d5}
<\  ,\ >_{L_\A,f,m}:\frac{ K_{m}}{  K_{m+1}} \times \frac{
K_{m}}{ K_{m+1}} \longrightarrow \R, \end{equation}
 whose signature
$\sigma_{\A,f,m}$ is independent of the linear map $L_\A$ chosen.
\end{Theorem}

\vskip 3mm
In
Section 2 we give a variational interpretation of  Theorem \ref{th1}.
Consider germs of analytic functions $f,f_1,f_2,\ldots,f_n$ in $\R^n$
such that $f,f_2,\ldots,f_n$ and $f_1,\ldots,f_n$
are regular sequences as holomorphic functions.
We consider the 1-parameter family of ideals
$
(f-t,f_2,\ldots,f_n).
$
Choose a small neighborhood $U_\C$  of $0 \in \C^{n}$  and a small $\varepsilon>0$
 such that:
\vskip 4mm
1) The sheaf of algebras on $U_\C$ defined by
$$
\BB_\C := \frac{{\mathcal O}_{U_\C}}{ (f_2,\ldots,f_{n})}
$$
is the structure sheaf of a 1-dimensional complete
intersection $\CC_\C \subset U_\C$ such that the map
\begin{equation}\label{finite}
f:\CC_\C \rightarrow \Delta_\varepsilon
\end{equation}
to the disk $\Delta_\varepsilon$ of radius $\varepsilon$ in $\C$
is a finite analytic map, the sheaf $f_*\BB_\C$ is a free
$\O_{\Delta_\varepsilon}$-sheaf of rank $\nu$ and $f^{-1}(0)=0$.

\vskip 2mm 2) $f_1|_{{\CC_\C}-\{0\}}$ is non-vanishing. \vskip 2mm
These conditions  can be fulfilled due to the
regular sequence hypothesis (\cite{E}, \cite{G}). Denoting by $f_*\BB^+$ the
sheaf on $(-\varepsilon,\varepsilon) \subset \R$ whose sections are the
fixed points of the conjugation map $\bar{\ }:f_*\BB_\C
\longrightarrow f_*\BB_\C$, we have that
 $f_*\BB^+$ is a free  $\AA_{(-\varepsilon,\varepsilon)}$-sheaf of rank
$\nu$. Its stalk over $0$ is
$$
\BB := (f_*\BB)_0^+  =
\frac{{\mathcal A}_{\R^n,0}}{ (f_2,\ldots,f_{n})} .
$$
Hence $\BB$
is a free ${\mathcal A}_{\R,0}$-module of rank $\nu$.
Introduce the 1-parameter family of $\R$-algebras
obtained by evaluation
\begin{equation}\label{3.3}
\B_{t_0}= f_*\BB^+ \otimes_\R \frac{\R[t]_{(t-t_0)}}{(t-t_0)}=
\left[\bigoplus_{p\in\ZZ_\C\cap
f^{-1}(t_0)}\frac{\O_{\C^n,p}}{(f-t_0,f_2,\ldots,f_n)}\right]^+.
\end{equation} $\B_0$ is a local algebra, $\B_{t_0}$ is a multilocal
algebra and they form a vector bundle of rank $\nu$ over
$(-\varepsilon,\varepsilon),$ whose sheaf of real analytic
sections is $f_*\BB^+$.

We define in the sheaf of sections $f_*\BB^+,$ a bilinear map
$$
<\ ,\ >:{f_*\BB^+} \times {f_*\BB^+}\
{\stackrel{\cdot}\longrightarrow }{f_*\BB^+}\ {\stackrel{\mathcal
L}\longrightarrow } {\mathcal A}_{(-\varepsilon,\varepsilon)}
\hskip 2mm , \hskip 2mm <a,b> = {\mathcal L}(a \cdot b),
$$
obtained by first applying the multiplication in
the sheaf of algebras $f_*\BB^+$ and then applying a chosen
${\mathcal A}_{(-\varepsilon,\varepsilon)}$-module map ${\mathcal
L}:{f_*\BB}^+ \longrightarrow  {\mathcal
A}_{(-\varepsilon,\varepsilon)}$
 having the property
that evaluating it at $0$ gives a  linear map $L_{\B_0}:\B_0
\longrightarrow \R$, verifying $L_{\B_0}([J]_{\B_0})>0$. The
evaluation of $<\ ,\ >$ at a fiber $\B_t$ is a bilinear form
defined on $\B_t$ and denoted by
$<\ ,>_t.$

 This family of non-degenerate bilinear forms is  the usual tool in the  Eisenbud--Levine
and Khimshiashvili theory (\cite{EL},\cite{K}) to calculate
the degree of the smooth map given by
$(f,f_2,\ldots,f_n):(\R^n,0)\to(\R^n,0)$.

\vskip 3mm
Define a sheaf map by multiplication with $f_1$
$$M_{f_1}:f_*\BB^+ \longrightarrow f_*\BB^+
\hskip 1cm M_{f_1}(b)=f_1b
$$
and a family of bilinear maps, that
we  call relative:
\begin{equation}\label{1.3c}
<\ ,\ >^{rel}:{f_*\BB}^+ \times {f_*\BB}^+\ \longrightarrow
{\mathcal A}_{(-\varepsilon,\varepsilon)} , \hskip 1cm <a,b>^{rel}
= <M_{f_1}(a),b>
\end{equation}
\begin{equation}\label{1.3d_t}
<\ ,>^{rel}_t:\B_t \times \B_t \longrightarrow \R ,
\hskip 1cm <[a]_t,[b]_t>^{rel}_t = <M_{[f_1]_t}([a]_t),[b]_t>_t.
\end{equation}
The bilinear forms $< , >^{rel}_t$  are
non-degenerate, for $t\neq0$, having signature $\tau_\pm$, for
$\pm t >0$. The form $<\ ,>^{rel}_t$ degenerates
for $t=0$ on $Ann_{\B_0}([f_1]_{\B_0})$ \cite{GGM}.
Expanding in Taylor series at $0$ the family of relative bilinear forms
we arrive at the setting in Jantzen \cite{J} and Vogan \cite{Vo}, where it is shown
how to obtain a flag of ideals
\begin{equation}
\label{flag2}
\ldots  \subset \tilde K_r \subset \ldots  \subset \tilde K_1  \subset \tilde K_0 = \B_0
\end{equation}
and bilinear forms in them and show how  to reconstruct from the signatures $\tau_m$ of these
bilinear forms the signatures $\tau_\pm$ (see Proposition \ref{Proposition 2.1}).
In our algebraic setting, the flag and the bilinear forms have the
algebraic description:

\vskip 3mm
\begin{Theorem}\label{Proposition 2.2} For the
family of bilinear forms $<\ ,\ >^{rel}_t,$ in the family of
algebras $\B_t$ (\ref{1.3d_t}) we have:
\begin{enumerate}
\item The set of $b \in \BB$ such that the function $t \rightarrow
<[b]_t,[b^\prime]_t>^{rel}_t$ vanishes at $0$ up to order $m$, for every $b^\prime
\in \BB$ is the quotient ideal
$$
(f^m:f_1) := \{b \in \BB \ / \ f_1 b \in (f^m) \} \subset \BB
$$
and
\begin{equation}
\label{flag}
\tilde{K}_m =
\frac{(f^m:f_1)}{(f) \cap (f^m:f_1)} \subset \frac{\BB}{(f)} =
\B_0.
\end{equation}
\vskip 2mm \item $(f) \cap (f^m:f_1) =
M_f((f^{m-1}:f_1)).$ \vskip 2mm \item  The bilinear form
$(b,b^\prime) \longrightarrow
L_{\B_0}([\frac{f_1b}{f^m}b^\prime]_{\B_0})$
\begin{equation}\label{z1} (f^m:f_1)\oplus (f^m:f_1)
{\stackrel{{\frac{f_1}{f^m}\cdot}} \longrightarrow} (f^m:f_1)
{\stackrel{ \tilde \pi_0} \longrightarrow} \frac{\BB}{(f)}=\B_0
{\stackrel{ L_{\B_0}} \longrightarrow }\R,
\end{equation}
where $\tilde \pi_0$ is the projection from $\BB$ to
$\BB/(f)=\B_0$, vanishes on $(f) \cap (f^m:f_1)$
and induces Jantzen's  bilinear form
\begin{equation}
\label{flag1} <\ ,\ >^{m}: \tilde K_{m} \otimes \tilde K_{m}
\longrightarrow \R \hskip 1cm <\ ,\ >^{m} =
<\frac{f_1\cdot}{f^m},\cdot>_{_0},
\end{equation}
giving the formula
$$ \tau_+ = \sum_{m \geq 0} \tau_m \hskip 1cm,\hskip 1cm \tau_- = \sum_{m \geq 0} (-1)^m \tau_m$$
\end{enumerate}
\end{Theorem}

\vskip 3mm
In Section 3, we show:

\vskip 3mm \begin{Theorem}\label{Proposition 3.1} There is an isomorphism
$\varphi:\tilde K_1 \rightarrow  K_1$, induced by multiplication
with the function $\frac{f_1}{f}$, which is sending the  flag
$\{\tilde K_m\}_{m\geq1}$ in $\B_0$ in (\ref{flag}) to the flag  $\{
K_m\}_{m\geq1}$ in $\A$ in (\ref{a}) and  Jantzen's bilinear forms
(\ref{flag1}) to the bilinear forms (\ref{d5}). Hence, for $m\geq
1,$ we have equal signatures $\tau_m = \sigma_{\A,f,m}$ and
$$ \tau_+ =  \tau_0 + \sum_{m  = 1}^{\ell+1} \sigma_{\A,f,m} \hskip 1cm,\hskip 1cm \tau_- = \tau_0 +\sum_{m =1}^{\ell +1} (-1)^m \sigma_{\A,f,m}.$$
\end{Theorem}

\vskip 3mm In Section 4 we apply these considerations for
calculating indices of vector fields. If $X=\sum_{i=1}^n
X^i\frac{\partial}{\partial x_i}$ is a real analytic vector field
with an algebraically isolated zero at $0$ in $\R^n$, then the
(Poincar\'e-Hopf) index of $X$ at $0$ is the signature of the
bilinear form $(\ref{d3})$ constructed for the finite dimensional
algebra
$$
\label{d11} \B:=\frac{{\mathcal A}_{\R^{n},0}} {
(X^1,\ldots,X^n)} \hskip 1cm,\hskip 1cm
<\  , \,>_{L_\B}: \B \times \B {\stackrel{\cdot}\rightarrow} \B
{\stackrel{L_\B }\rightarrow} \R
$$
where $L_\B:\B \to \R$ is a linear map with $L_\B(J_X)>0$ (see
\cite{EL} and \cite{K}). Now assume further that $f:(\R^n,0) \to
(\R,0)$ is a real analytic function, that $X$ is tangent to the
fiber $V_0 := f^{-1}(0)$, giving the relation $df(X)=hf$ with $h$ a
real analytic function called the cofactor. If $0$ is a smooth point
of $V_0$ then the signature $\sigma_{\B,h,0}$ of the order $0$
bilinear form
$$<\ ,\ >_{L,h,0}:\frac{\B}{Ann_{\B}(h)} \times
\frac{\B_\R}{Ann_{\B}(h)}
{\stackrel{\cdot}\longrightarrow}
\frac{\B}{Ann_{\B}(h)}
{\stackrel{L}\longrightarrow}
\R
$$
$$
L:\frac{\B_\R}{Ann_{\B}(h)}
\longrightarrow \R \hskip 3mm , \hskip 3mm L(\frac{J_\B}{h})>0
$$
is the Poincar\'e-Hopf index at $0$ of the vector field $X|_{V_0}$,
as can easily be deduced using the implicit function theorem. If $0$
is an isolated critical point of $V_0$ and the dimension $n$ of the
ambient space is even, in \cite{GM2} it is proved that
\begin{equation}\label{indexeven}
Ind_{V_{+,0}}(X)=Ind_{V_{-,0}}(X)=\sigma_{\B,h,0}-\sigma_{\A,h,0}.\end{equation}
If $n$ is odd, it is proved in \cite{GM1} that

\begin{equation}\label{indexodd}
Ind_{V_{\pm,0}}(X)=\sigma_{\B,h,0}+K_{\pm}.\end{equation}
 In  the case of odd dimensional ambient space and for $f$ a germ of a real
analytic function with an algebraically isolated singularity at $0$,
we calculate the constants $K_\pm$ by studying the family of contact
vector fields

$$
 X_t= (f-t)\frac{\partial}{\partial x_1}+ \sum_{i=1}^N[
\frac{\partial f}{\partial x_{2i+1}}\frac{\partial }{\partial
x_{2i}}- \frac{\partial f}{\partial x_{2i}}\frac{\partial}{\partial
x_{2i+1}}],
$$
where $f_j := \frac{\partial f}{\partial x_j}$. For $t \neq 0$, the
signatures of the relative bilinear forms correspond to the sum of
the Poincar\'e-Hopf indices of the restriction of $X_t$ to $V_t$.
Our transport from the algebra $\B_0$ to the Jacobian algebra $\A$
is a local analogue of the Poincar\'e-Hopf Theorem relating
information of the singular point of $X$ to invariants of the
singularity of $f$.

\vskip 3mm
Using these explicit computations for contact vector
fields, we conclude the search for an algebraic formula for the
real GSV-index using local algebra by determining the values of the constants
$K_\pm$:

\vskip 0.5cm \begin{Theorem}\label{th2} Let $V$ be an algebraically
isolated hypersurface singularity in $\R^{2N+1}$, then the constants
$K_\pm$ in    (\ref{indexodd}) relating the GSV-index and the
signature $\sigma_{\B,h,0}$ are:
$$
K_+ = \sum_{m \geq 1} \sigma_{\A,f,m} \hskip 1cm,\hskip 1cm
K_-  =  \sum_{m \geq 1} (-1)^m\sigma_{\A,f,m}.$$
\end{Theorem}

\vskip 5mm

\section{Higher Order Signatures in $\A$}

We use the definitions and notations of the Introduction. The
algebra $\A$ has an intrinsic $\A$-valued bilinear map, which is
the multiplication in $\A$:
\begin{equation}\label{a*}
(\ ,\ )_\A:\A \times \A \longrightarrow \A \hskip 2cm  (a ,b
)_\A:=ab.
\end{equation}
In terms of this pairing and the non-singular pairing
$<\ ,>_{L_\A}$ in (\ref{d3}), the orthogonal of an ideal $I \subset \A$
is the annihilator ideal in the algebra $\A$:$\  I^\perp = Ann_\A(I)$.
The process of taking the orthogonal induces an involution in the
set of ideals of $\A$, which is reversing the natural inclusions
of sets in $\A$. In particular, the orthogonal to the maximal
ideal is the socle. It is 1-dimensional and the class of the
Jacobian $J_\A$  is a generator (see \cite{EL}, \cite{K},
\cite{E}).
\vskip 3mm
Choose now an element $f \in \A$ in the
maximal ideal. Consider the linear map induced in $\A$ by
multiplication with $f$:
$$
M_f:\A \longrightarrow \A \hskip 2cm M_f(a)=f a.
$$
For $j \geq 1$, the maps $M_f^j$ are  selfadjoint maps for the
bilinear map (\ref{a*}):
$$
(M_f^ja,b)_\A = f^jab = af^jb = (a,M_f^jb)_\A ,
$$
and hence they are also  selfadjoint maps for the bilinear form
$<\ ,\ >_{L_\A}$. We have that
$$
Ann_\A(f^j) = Ker(M_f^j)\hskip 1cm  \hbox{and} \hskip1cm (f^j) = Im(M_f^j)
$$
and each of these spaces is the  orthogonal of the other in $\A$ ,
since $M_f^j$ is selfadjoint.

\vskip 3mm Consider the  flag of ideals in $\A$
\begin{equation}\label{1.1}
0 \subset (f^{\ell}) \subset (f^{\ell-1}) \subset \cdots \subset
(f^2) \subset (f) \subset \A,
\end{equation}
where $\ell$ is
minimal with $f^{\ell+1}=0$ and the orthogonal flag of ideals
\begin{equation}\label{1.2} 0 \subset Ann_\A(f) \subset Ann_\A(f^2)
 \subset \cdots \subset Ann_\A(f^{\ell-1}) \subset
Ann_\A(f^{\ell})\subset \A. \end{equation}
The linear map $M_f:\A
\longrightarrow \A$ is a nilpotent map $M_f^{\ell+1}=0$.

\vskip 3mm
\begin{Lemma}\label{Lemma 1.1} For $j=1,\ldots,\ell +1,$ there are linear subspaces $P_j$
of $\A$, called primitive subspaces, such that
\begin{equation}\label{aa}
\A=\bigoplus_{j=1}^{\ell +1} [\bigoplus_{k=0}^{j-1}M_f^kP_j],
\end{equation}
with $M_f^{j-1}:P_j \longrightarrow \A$ injective and
$M_f^{j}(P_j)=0$. The mapping $M_f:\A\to\A$ is in Jordan
canonical form in any basis obtained by choosing bases of each of
the spaces $P_j$ and extending them to a basis of $\A$ by the
action of $M_f$ as in (\ref{aa}).
\end{Lemma}
{\bf Proof:} We recall how to choose a basis of $\A$ as a vector
space over $\R$ that expresses $M_f$ in Jordan canonical form.
Inductively, let us begin by choosing linearly independent vectors
$v_1, \ldots, v_{n_{\ell +1}}$ generating a vector space $P_{\ell +1}$
complementary to $Ann_\A(f^{\ell})$ in $\A$ and choose as
first vectors of a basis of $\A$ the vectors
$$
\{v_j,fv_j,\ldots,f^{\ell}v_j\}_{j=1,\ldots,n_{\ell+1}}.
$$
With $P_{\ell +1}$ we construct the Jordan blocks of maximal size
$\ell$ of $M_f$. Then, we choose linearly independent vectors
$v_{n_{\ell +1}+1},\ldots,v_{n_{\ell+1}+n_{\ell}}$ generating a
vector space $P_{\ell}$ with the property that
$$
Ann_\A(f^{\ell-1}) \oplus M_f( P_{\ell +1}) \oplus P_{\ell} = Ann_\A(f^{\ell}).
$$
We choose the next part of the basis by choosing the vectors
$$
\{v_j,fv_j,\ldots,f^{\ell-1}v_j\}_{j=n_{\ell+1}+1,\ldots,n_{\ell+1}+n_{\ell}}
$$
to construct the Jordan blocks of size $\ell-1 $, and so on. The
space of $1$-st primitive vectors $P_1$ is formed  of vectors in
$\A$ with the property that
$$
M_f^\ell(P_{\ell+1}) \oplus M_f^{\ell-1}(P_{\ell}) \oplus
\cdots\oplus
 M_f^2(P_3)\oplus
M_f(P_2)\oplus P_1 = Ann_\A(f).
$$
\qed

\vskip 2mm We call the vectors in $P_j$ $j$th-primitive vectors,
and we denote by $n_j$  the dimension of $P_j$. Hence $n_j$ is
also the number of Jordan blocks of size $j$ in $M_f$. It is
convenient to present the direct sum decomposition (\ref{aa}) by
the matrix:
\begin{equation}\label{b}
\A = \left(\begin{matrix} P_1 & P_2 & P_3 & P_4 & \cdots &
P_{\ell} & P_{\ell+1} \cr 0 & M_fP_2 & M_fP_3 & M_fP_4 & \cdots &
M_fP_{\ell} & M_fP_{\ell+1} \cr 0 & 0 & M_f^2P_3 & M_f^2P_4 & \cdots
& M_f^2P_{\ell} & M_f^2P_{\ell+1} \cr
 &  &  &  & \cdots &  & \cr
0 & 0 & 0 & 0 & \cdots & M_f^{\ell-1}P_{\ell} &
M_f^{\ell-1}P_{\ell+1} \cr 0 & 0 & 0 & 0 & \cdots & 0 & M_f^\ell
P_{\ell+1} \cr,
\end{matrix}\right)
\end{equation}
 meaning that an element of $\A$ has components in the form of
an upper triangular matrix where the $(i,j)^{th}$-entry of the
matrix is an arbitrary element in $M_f^{i-1}(P_{j})$, with
$i, j=1,\ldots,\ell+1$. Each column is formed by equidimensional
subspaces, until we reach the zero subspace,
and the map $M_f$ acts as a map preserving columns and
descending one row. Hence, restricting to a column in (\ref{b}),
the map $M_f$ is an isomorphism until it reaches the diagonal,
where $M_f$ is the zero map. \vskip 3mm Using this representation
for $\A$ and recalling the flag of ideals (\ref{a}), we have:

\vskip 3mm

\begin{Lemma}\label{Lemma 1.2} 1) The ideal $(f^m)$
 is formed by the last $\ell+1-m$ rows of the matrix
(\ref{b}).
\vskip 3mm
2)  Its orthogonal $Ann_\A(f^m)$ is formed
by the elements in a band of width $m$ above the diagonal in
(\ref{b}), including the diagonal.
\vskip 3mm
3) The ideal $ K_m$
in (\ref{a}) is formed by  the lower $\ell+2-m$ diagonal terms.
\vskip 3mm
4)  The  ideal $ K_m^\perp$, orthogonal to $  K_m$ is
\begin{equation}\label{b1}
K^\perp_m = (f) + Ann_\A(f^{m-1}).
\end{equation}
\end{Lemma}
\begin{Example}\label{example1.1}
For $\ell=3$ and $m=3$, we have
$$
(f^2) = \left(\begin{matrix} 0 & 0 & 0 & 0 \cr 0 & 0 & 0 & 0 \cr 0
& 0 & M_f^2P_3 & M_f^2P_4 \cr 0 & 0 & 0 & M_f^3P_4 \cr
\end{matrix}\right)
\hskip 1cm Ann_\A(f^2) = \left(\begin{matrix} P_1 & P_2 & 0 & 0
\cr 0 & M_fP_2 & M_fP_3 & 0 \cr 0 & 0 & M_f^2P_3 & M_f^2P_4 \cr 0
& 0 & 0 & M_f^3P_4 \cr
\end{matrix}\right)
$$
$$
 K_3 = \left(\begin{matrix} 0 & 0 & 0 & 0 \cr 0 & 0 & 0 & 0 \cr 0 & 0 &
M_f^2P_3 & 0 \cr 0 & 0 & 0 & M_f^3P_4 \cr
\end{matrix}\right)
\ \ ;\ \ (f) + Ann_\A(f^2) = \left(\begin{matrix} P_1 & P_2 & 0 &
0 \cr 0 & M_fP_2 & M_fP_3 & M_fP_4 \cr 0 & 0 & M_f^2P_3 & M_f^2P_4
\cr 0 & 0 & 0 & M_f^3P_4 \cr
\end{matrix}\right) =  K_3^\perp.
$$
\end{Example}
\vskip 3mm
\noindent{\bf Proof of Lemma \ref{Lemma 1.2}}. Since
$M_f$ corresponds to going down 1 row in (\ref{b}), parts 1, 2 and
3, are clear. To prove part 4, note first that
 $$ (f) + Ann_\A(f^{m-1}) \subset  K_m^\perp.$$
The ideal $(f)$ is given by all the terms in (\ref{b}), except for
the first row. Since $Ann_\A(f^{m-1})$ is the band matrix above the
diagonal of width $m-1$, we obtain that the only contribution of
$(f) + Ann_\A(f^{m-1})$ to $(f)$ is given by the first $m-1$ terms in
the first row. On the other hand ${K}_m=Ann_\A(f) \cap (f^{m-1})$
consist of the last $\ell+2-m$ terms in the diagonal. We observe
on using (\ref{b}) that the ideals $(f) + Ann_\A(f^{m-1})$ and $K_m$
have complementary dimensions in $\A$. Now (\ref{b1}) must hold,
as the bilinear form
$<\ ,>_{L_\A}$ is non-degenerate. \qed

\vskip 3mm \begin{Proposition}\label{Proposition 1.1} For the
bilinear forms in (\ref{aaa}), we have: \vskip 2mm
\begin{enumerate}
\item{}  $<\ ,\ >_{{L_\A},f,0} = <f\cdot,\cdot>_{L_\A}$ has $ K_1
= Ann_\A(f)$ as
 degeneracy locus and the induced non-degenerate bilinear form in
$\frac {\A }{Ann_\A(f)}$ is obtained by choosing $\frac{J_\A}{ f}$
as generator of the 1 dimensional socle of $\frac{\A}{ Ann_\A(f)}$
and defining the bilinear form as multiplication followed by a
real valued map sending $\frac{J_\A}{ f}$ to a positive number.
\vskip 2mm
\item{} The bilinear form $<\ ,\ >_{L_\A,f,1} = <\ ,\
>_{L_\A}|_{K_1 \times K_1}$ has $ K_2$
  as degeneracy locus.
\vskip 2mm \item {} For $m \geq 2$ the bilinear form
$<\ ,>_{L_\A,f,m}$ in (\ref{aaa})
is well defined and has $ K_{m+1}$
 as degeneracy locus.
\end{enumerate}
\end{Proposition}
\noindent{\bf Proof:}  1)  The inner product
$<f\cdot,\cdot>_{L_\A}$ vanishes on $Ann_\A(f)$. If
$<fa,a^\prime>_{L_A}=0$ for all $a^\prime$, then $fa=0$, since
(\ref{d3}) is a non-degenerate bilinear form on $\A$ (\cite{E},
\cite{EL}, \cite{K}). Hence, $<\ ,\ >_{L_\A,f,0}$ has $K_1$ as
degeneracy locus. The algebra $\frac{\A}{Ann_\A(f)}$ has a
one-dimensional socle generated by the class of $J_\A/f$ (see
\cite{GGM}  for more details).
\vskip 3mm
2)  Note first
that $K_2=Ann_\A(f) \cap (f)=K_1 \cap K_1^\perp$, by (\ref{a}) and
(\ref{b1}). Hence, given $a\in K_2$ and any $b\in K_1$, it follows
that $(a,b)_{\A}=0,$ so $K_2$ is contained in the degeneracy locus
of $<\ ,\ >_{L,f,1}$. On the other hand, let $a \in  K_1 - K_2=
K_1 -  K_1^\perp$. Then $a K_1$ is a non-zero ideal in $\A$, and so
contains the socle of $\A$. We obtain an expression $J_\A=ac$, for
some $c \in K_1$. Hence, $<a,c>_{L_\A}=L_\A(ac)=L_\A(J_\A)>0$, so
that $a$ is not in the degeneracy locus of $<\ ,\ >_{L,f,1}$.
\vskip 3mm
3)  Let $m \geq 2$. We first show that the bilinear
form $<\ ,>_{L_\A,f,m}$ is well defined, i.e. is independent of the division
by $f^{m-1}$ in $K_{m}$. Let $a$, $b$ be in
$K_{m}=Ann_\A(f)\cap(f^{m-1})$. Then there exists $a_1\in \A$
such that $a=a_1f^{m-1}$ and $<a,b>_{L_\A,f,m}= \hfill  <a_1,b>_{L_\A}$.
Let also $a=a_2f^{m-1}$. Then $<a_1,b>_{L_\A}=\quad <a_2,b>_{L_\A}$,
because $a_1-a_2\in Ann_{\A}{(f^{m-1})}$ and $b\in(f^{m-1})$.
\vskip3mm
If $a \in  K_{m+1}$, then $\frac{a}{f^{m-1}} \in (f)$, and since $b \in
K_{m} \subset Ann_\A(f)$, we have $\frac{a}{ f^{m-1}}b=0$. Hence, the form
$<\, \ >_{L_\A,f,m}$
degenerates on $ K_{m+1}$.
\vskip 3mm
Let $a  \in  K_{m} - K_{m+1}.$ In order to prove that the form $<\ , \ >_{L_\A,f,m}$ is
non-degenerate on $a$, we have to show that $\frac {a}{ f^{m-1}}\not\in K_{m}^\perp$.
 Using the representation (\ref{b}), and part 3) of Lemma~\ref{Lemma 1.2},
the $a_{m,m}$ entry  in $a$ is  not zero, and $a_{m,m}\in
M_f^{m-1}P_{m}$.
 Now $\frac{a}{ f^{m-1}}$ is obtained by lifting all the
elements in the representation by $m-1$ rows, keeping the columns
fixed. We observe that $\frac{a}{ f^{m-1}}\not\in(f)$. It now
suffices to show that $\frac{a}{ f^{m-1}}\not\in Ann_\A(f^{m-1})$.
But by part 4) of Lemma~\ref{Lemma 1.2}, the space
$Ann_\A(f^{m-1})$ is given by the band matrix of width $m-1$, including the
diagonal.
Hence, $\frac{a}{ f^{m-1}}$ is not an element of
$Ann_\A(f^{m-1}).$
\qed

\vskip 3mm
\noindent{\bf Proof of Theorem \ref{th1}}. By
Proposition~\ref{Proposition 1.1}, we have that $K_{m+1}$ is the
 locus of the bilinear form $<\ ,>_{L_\A,f,m}$, so that
$K_{m}/K_{m+1}$ inherits a non-degenerate bilinear form. The linear
forms $L_\A$, verifying $L_\A(J_\A)>0$ form an open connected set in
the dual space $\R^{n^*}$. The signature is an integer valued
continuous function of $L_\A$, hence it is constant. \hfill $\qed$

\vskip 3mm \begin{Corollary}\label{Corollary 1} For $m\geq1$, the
mapping \begin{equation}\label{iso} M^{m-1}_f:P_{m}
\longrightarrow K_{m}/K_{m+1}\end{equation}
 is a well defined
isomorphism.
 The
pairing of $m$-primitive vectors
$$
<\ ,\ >^{prim}_{L_\A,m}:P_{m} \times P_{m} \longrightarrow \R,
\hskip 2cm <a,b>^{prim}_{L_\A,m}:= <M_f^{m-1}a,b>_{L_\A}
$$
 is a non-degenerate symmetric bilinear pairing,
induced by the pairing (\ref{d5}) via the isomorphism (\ref{iso}).
\end{Corollary}
\noindent{\bf Proof:} Using the representation (\ref{b}) for the
elements of $\A$ and the description of $K_m$ given in part 3 of
Lemma~\ref{Lemma 1.2}, we have that $M^{m-1}_f:P_{m}
\longrightarrow K_{m}$ is injective, and hence (\ref{iso}) is a
well defined isomorphism. The pull back of the non-degenerate
bilinear form on $K_{m}/K_{m+1}$ via this last map is
$$
<b,b'>_{L_\A,m}=L_\A(f^{m-1}\cdot b\cdot
b')=L_\A((\frac{1}{f^{m-1}}\cdot f^{m-1}b)\cdot f^{m-1} b')=$$
$$=<M^{m-1}_f(b),M^{m-1}_f(b')>_{L_\A,f,m}.$$ \qed
\begin{Example}\label{Example 1.2} Let $f$ be a germ of
a real analytic function in $\R^n$, with an algebraically isolated
critical point. This means that the ideal generated by the partial
derivatives of $f$ in the ring of germs of holomorphic functions
has finite codimension.
 Let $\A=\A(f)$   be given by (\ref{d1}), with $f_i := \frac{\partial f}{\partial
 x_i}$ and let $\ell$ be as in (\ref{1.1}).
Let $\sigma_{f,m}=\sigma_{\A(f),f,m},$ $m=0,\ldots,\ell+1$, be the
signatures given by Theorem~\ref{th1}. These are invariants
associated  to the germ $f$. We call $\sigma_{f,m}$ the order $m$
signature of $f$.
\end{Example}

\vskip 1cm
\section{The Family of Bilinear Forms in $\BB$.}
In this section we construct a family of bilinear forms
$<\ ,>^{rel}$, that we call relative, which is constructed from the
equations
$$f-t=f_2=\cdots=f_n=0$$
which are non-degenerate for $t \neq 0$. We do Taylor series
expansion of $<\ ,\ >^{rel}$ and determine an algebraic procedure
to compute the signatures for $t\neq 0$ in terms of local linear
algebra in the ring $\frac {{\mathcal A}_{\R^{n},0}} {
(f_2,\ldots,f_{n})}$ via the first terms of the above Taylor
series expansion. \vskip 5mm

Recall the setting and definitions of the Introduction.
Note in particular that
since the map in (\ref{finite}) is a finite analytic map, the inverse image $f^{-1}((-\varepsilon,\varepsilon))$
is a finite union of  curves (parameterized by $(-\varepsilon,0]$ or $[0,\varepsilon)$),
which come together at $0$. The conjugation map permutes them, and the
fixed components correspond to  $\CC :=\CC_\C \cap \R^{n}$.
Hence $\CC$ consists either of $0$ only or of a finite number of these real curves
all passing through $0$, which is its only singular point. Note that the degree of the covering map
$f:\CC - \{0\} \longrightarrow (-\varepsilon,\varepsilon)-\{0\}$ may be distinct for $t>0$ and $t<0$.
In the sheaf $f_*\BB^+$
we have information about the points $\{f=t\}_{t \in (-\varepsilon,\varepsilon)} \cap \ZZ_\C$ in $U_\C$, real or complex.

\vskip 5mm
\begin{Lemma}\label{Lemma 2.1} The signature of the
non-degenerate bilinear forms $<\ ,\ >_t$ on $\B_t$ is independent
of $t$ and it is equal to the sum of the signatures of the
bilinear forms computed on the local rings $\B_{t,p}$ for $p \in
\ZZ\cap f^{-1}(t)$, for each $t \in (-\varepsilon,\varepsilon)$.
\end{Lemma}
\noindent{\bf Proof:} This is the usual procedure due to Eisenbud--Levine
and Khimshiashvili (\cite{EL},\cite{K}) to calculate
the degree applied to the smooth map given by
$(f,f_2,\ldots,f_n):(\R^n,0)\to(\R^n,0)$. In particular, the
contribution to the signature coming from points in $\ZZ_\C - \ZZ$
is always 0. \qed

\vskip 3mm

Using a trivialization of $f_*\BB^+$,  we can transfer the
relative forms
$<\ , \,>^{rel}_t$ from $\B_t$ to $\B_0$. So we have
a family of bilinear forms
that we denote by $<\ ,\ >_t$.
We are interested in the Taylor expansion of this family of
bilinear forms at $t=0$.  We will use the following
Proposition, containing results from  Jantzen \cite{J} and Vogan \cite{Vo}:
\begin{Proposition}\label{Proposition 2.1}
Let $< \,, >_t$,
$t\in(-\varepsilon,\varepsilon)$, be an analytic family of
 forms on a finite dimensional vector space $\B_0$. Assume
that the forms $< \ ,\ >_t$ are nondegenerate for $t\ne0$. Let
$\tilde K_i$, $i=0,\ldots,r$, be the set of $[b]_0 \in\B_0$,
such that the functions $t\mapsto<[b]_0,[b^\prime]_0>_t$
vanish at $0$ up to order $i$ for any $b^\prime \in \BB$. Then \vskip 2mm
\begin{enumerate}
\item{} For $i=0,\ldots,r,$ a bilinear form $< \ ,\ >^i$,  is well
defined on $\tilde K_i$ by
\begin{equation}\label{bbb}
<[b]_0,[b^\prime]_0>^i= \frac{1}{i!}\frac{d^i}{dt^i}<b,b^\prime>_t|_{t=0}
.\end{equation}
\vskip 2mm
\item{} The bilinear form $<\ ,\ >^i$
degenerates on $\tilde K_{i+1}$, and induces a  nondegenerate
bilinear form on $\tilde K_i/\tilde K_{i+1}$. Denote  its
signature by $\tau_i$. \vskip 2mm
\item{}
The signatures $\tau_+$ and $\tau_-$ of the forms
$<\ ,\ >_t$ on $\B_0$, for $t>0$ and $t<0$, respectively, are given by

\begin{equation}\label{z0}
\tau_+=\sum_{i=0}^r\tau_i,\hskip 1cm \hskip 1cm
\tau_-=\sum_{i=0}^r(-1)^i\tau_i. \end{equation}
\end{enumerate}
\end{Proposition}
\vskip 3mm

\noindent{\bf Proof of Theorem \ref{Proposition 2.2}:}
 1)  Let $b \in \BB-\{0\}$. Then, there is a unique integer $j$
and $c \in \BB$ with $[c]_{\B_0}\neq0$ such that $f_1b = f^jc$.
Since $[c]_{\B_0}\neq0,$ and $[J]_{\B_0}$ is a generator of the
socle of $\B_0$ we may find $e_0 \in \B_0$ such that
$[c]_{\B_0}e_0 = [J]_{\B_0}$. Choose any $e \in \BB$ with the
property that $[e]_{\B_0}= e_0$, so that $\L(ce)(0)=
L_{\B_0}([c]_{\B_0}e_0) =  L_{\B_0}([J]_{\B_0})
 \neq0$. For any $b^\prime \in \BB$, we have
$f_1bb^\prime =f^jcb^\prime \in (f^j)$. Hence,
$$<b,b^\prime>^{rel}= \L(f_1bb^\prime)=\L(f^jcb^\prime)=t^j\L(cb^\prime) \in (t^j)$$
and
$$<b,e>^{rel}= \L(f_1be)=\L(f^jce)=t^j\L(ce) \in (t^j)-(t^{j+1}).$$
Hence, if $b$ is as in the statement of part 1), we have that
$j\geq m$ and hence $f_1b \in (f^m)$, i.e. $b \in (f^m:f_1)$. This
proves the first assertion. The second assertion follows from the
first by evaluating it at $t=0$ and using (\ref{3.3}).

\vskip 2mm
2) Let $b \in (f) \cap (f^m:f_1)$. Then $b=cf$ and
$(cf)f_1=ef^m$, so that $cf_1=ef^{m-1}$. Hence, $c \in
(f^{m-1}:f_1)$ and $b =cf \in M_f (f^{m-1}:f_1)$. The converse is
obvious.

\vskip 3mm
3) Let $b \in (f) \cap (f^m:f_1)$ and
$b^\prime \in (f^m:f_1)$. Then
$$\left(\frac{f_1b}{f^m}\right)b^\prime = b \left(\frac{f_1b^\prime}{f^m}\right) \in (f),$$
since $b \in (f)$. Hence $[<
{\frac{f_1b}{f^m}},b^\prime>]_{\B_0}=0$ and the bilinear form in
(\ref{z1}) vanishes on $(f) \cap (f^m:f_1)$. Taking the quotient
by $(f) \cap (f^m:f_1)$, we obtain by part 1) that it is a
bilinear form defined on $\tilde K_m$ and it has the same
expression as Jantzen's form, since $f^m=t^m$, so they coincide.
 $\qed$

\vskip 5mm \section {Transporting the signatures to the algebra
$\A$} The aim of this section is to establish a relationship
between the higher order bilinear forms $< , >_{L_A,f,m}$
(\ref{aaa}) and their signatures $\sigma_{\A,f,m}$ in the algebra
$\A$ and  Jantzen's relative forms
$< ,>^m$ (\ref{bbb}) and their signatures $\tau_m$ in $\B_0$.

Define the isomorphism of $\BB$-modules
$$\begin{matrix} \BB & & \BB \cr
\cup&&\cup\cr (f:f_1) & {\stackrel{{\Phi}} \longrightarrow} &
(f_1:f) \cr
\end{matrix}$$
$$\Phi( b) = \frac{ b f_1}{f}, \hskip 1cm  \hskip 1cm
\Phi^{-1}( c) = \frac{ c f}{f_1}$$
\begin{Lemma}\label{Lemma 3.1}  The isomorphism $\Phi$ induces
isomorphisms of $\BB$-modules,  for $m\geq 1$:
\begin{equation}\label{*1}\Phi: (f^m:f_1) \longrightarrow (f_1:f) \cap
(f^{m-1})
\end{equation}
$$\Phi: (f^m) \longrightarrow  (f^{m-1}f_1)$$
$$\varphi: \tilde K_1=Ann_{\B_0}(f_1) \longrightarrow K_1 = Ann_\A(f).$$
\end{Lemma}

\noindent{\bf Proof:} If $ b \in (f^m:f_1)$, then there exists $
c\in\BB$ such that $ b f_1 =  c f^m$. Hence
$$\Phi( b) = \frac{ b f_1}{f}=  c f^{m-1} \in
(f_1:f) \cap (f^{m-1}).$$ Conversely, if $ c =  d f^{m-1} \in
(f_1:f),$ then
$$ d f^m =  c f =  e f_1 \hskip 1cm \Rightarrow \hskip 1cm
 e = \Phi^{-1}( c) \in (f^m:f_1).$$ This proves the
first assertion. The second one is just $\Phi( b f^m)= b
f^{m-1}f_1$. The third assertion is obtained by taking the
quotient of the first assertion in the Lemma by the second
relation in the case $m=1$ . \qed

\vskip 5mm

Let $f_2,\ldots,f_n$ be a regular sequence of holomorphic functions, denote
the volume form by $dVol = dx_1 \wedge \cdots \wedge dx_n$, and let $\ZZ_\C$ be the complete intersection
$f_2=\ldots=f_n$ as in Section 2. For any holomorphic function $g$ define the Jacobian of $g$ by
$$
dg \wedge df_2 \wedge \ldots \wedge df_n := Jac(g) \, dVol \ \ , \ \
Jac(g) =
\left| \begin{array}{ccc}
        \frac{\partial g}{\partial x_1} & \cdots & \frac{\partial g}{\partial x_n} \\
        \frac{\partial f_2}{\partial x_1} & \cdots & \frac{\partial f_2}{\partial x_n}\\
                            \cdot                        & \cdots & \cdot \\
        \frac{\partial f_n}{\partial x_1}  & \cdots & \frac{\partial f_n}{\partial x_n}
       \end{array} \right|.
$$
Recall the construction of the generator of the Rosenlicht differentials (see \cite{MvS}) or dualizing
module, which is a rational differential form $\omega_0$ on $\ZZ_\C$ having the property
$$\omega_0 \wedge (df_2\wedge\ldots\wedge df_n)|_{\ZZ_\C} = dVol|_{\ZZ_\C} \in
\frac{\Omega^n_{\C^n}}{(f_2,\ldots f_n)\Omega^n_{\C^n}}. $$
The dualizing module on $\ZZ_\C$ is then $\O_{\ZZ_\C} \omega_0$, and it consists of all
rational differential  forms $\sigma$ on $\ZZ_\C$ that have the property that the residue at $0$
of $h \sigma $ is $0$, for any holomorphic function $h$ on $\ZZ_\C$. Recall also that the residue
of a differential form $\sigma$ at $0 \in \ZZ_\C$ is the sum of the residues of the rational differential form
$\nu^*\sigma$ at $\nu^{-1}(0)$, where $\nu$ is the normalization map of $\ZZ_\C$.
Directly from the definitions above, one obtains that for any holomorphic function $g$ on $\C^n$
$$ d (g|_{\ZZ_\C}) = Jac(g)|_{\ZZ_\C} \omega_0.$$
Note that
  the logarithmic derivative of $g|_{\ZZ_\C}$ is $\frac{Jac(g)}{g}\omega_0$ and its
residue at $0$ is the sum of the vanishing orders
of the function $g \circ \nu$ at $\nu^{-1}(0)$, and hence a positive integer.

\vskip 3mm

\begin{Lemma}\label{Lemma 3.2}
Let $f,f_1,\ldots,f_n$, and $\varphi$ be as in Section 2. Let $J_{\B_0}$
and $J_\A$ be the Jacobians of $(f,f_2,\ldots,f_n)$ and
$(f_1,\dots,f_n)$ respectively.
   Then there exists
a positive constant $c=c(f)$ such that $\varphi(J_{\B_0})=c J_\A.$

\end{Lemma}

\noindent{\bf Proof: }
Since $([f]_\A) \varsubsetneq \A$, then taking
orthogonal of this relation, we obtain that the ideal $ K_1$ is
not the $0$-ideal, and hence $ K_1$ contains the socle.
 Since $\varphi:\tilde K_1 \longrightarrow
K_1$ is an isomorphism of non-zero ideals,
each containing its corresponding 1-dimensional
socle,  the map $\varphi$ sends the  socle
ideal to the corresponding  socle ideal. Hence
$\varphi$ sends the Jacobian of $\B_0$ to a non-zero multiple of
the corresponding Jacobian of $\A$.

\vskip 3mm

  Thus we know that there is a non-zero real number $c$ with the property
$$
\left[ \frac{f_1 Jac(f) }{f}\right]_\A  = c[Jac(f_1)]_\A.
$$
Hence there is a holomorphic function $h$ on $\ZZ_\C$ with the property
$$
\frac{f_1 Jac(f)}{f} - c Jac(f_1) = h f_1
$$
Dividing by $f_1$ and multiplying by $\omega_0$ we obtain
$$
\frac{Jac(f)}{f} \omega_0 - c \frac{Jac(f_1)}{f_1} \omega_0  = h \omega_0
$$
Taking residues at $0$ we obtain that
$$
n_1+\cdots+n_r -c (m_1+\cdots+m_r) = 0
$$
where the $n_i$ and $m_j$ are the vanishing orders of the functions $f \circ \nu$ and $f_1 \circ \nu$
at $\nu^{-1}(0)$, respectively.
Hence $c$ is a positive rational number. \qed

\vskip 3mm
The real valued bilinear forms on $\A$ and on $\B_0$ depended on
the choice of real valued linear functions $L_\A:\A
\longrightarrow \R$ and $L_{\B_0}:\B_0 \longrightarrow \R$ which
have the property of sending the
 corresponding Jacobians
to a positive number. Having chosen $\L$ and hence $L_{\B_0}$, we
will choose $L_\A$ subject to the compatibility condition
\begin{equation}\label{3.3.1}
L_{\B_0}|_{\tilde K_1}= L_\A \circ \varphi.\end{equation}

\vskip 5mm \noindent{\bf Proof of Theorem \ref{Proposition 3.1}:} Let $m
\geq 1$ and consider the commutative diagram:

$$
\begin{matrix} (f^m:f_1) \oplus (f^m:f_1) & {\stackrel{ \frac{f_1
\cdot\ }{f^m}\cdot} \longrightarrow} & (f^m:f_1) & \stackrel{
\tilde \pi_0 } \longrightarrow & \tilde K_1 & {\stackrel{ L_{B_0}
} \longrightarrow} & \R \cr \Phi\oplus \Phi\downarrow  & &
\downarrow \Phi & &\downarrow \varphi & & \ \ \downarrow Id \cr
(f_1:f)\cap(f^{m-1}) \oplus (f_1:f)\cap(f^{m-1}) & {\stackrel{
\frac{1\cdot\ }{f^{m-1}} \cdot \ } \longrightarrow} &
(f_1:f)\cap(f^{m-1}) & \stackrel{ \pi_0 } \longrightarrow & K_1 &
{\stackrel{ L_{\A} } \longrightarrow} & \R \cr
\end{matrix}.
$$
Here the mapping $\frac{f_1 \cdot }{f^m}\cdot$ acts on a couple
$$(a,b)\in (f^m:f_1) \oplus (f^m:f_1) \ \ \ \hbox{ by}
\ \ \ \ (a,b) \longrightarrow \frac{af_1  }{f^m}b,$$ and similarly
for $\frac{1\cdot }{f^{m-1}}\cdot$. The mapping $\pi_0$ is
obtained by reducing mod $(f_1)$. The vertical maps are
isomorphisms, so we may interpret the commutative diagram as
providing a conjugation of the top bilinear form into the bottom
bilinear form. We reduce the first row by $(f)$ and the second row
by $(f_1)$. This is possible since $\Phi(f)=f_1$ and both bilinear
forms degenerate in the submodules in the denominator of the
quotient. We thus obtain that the $m^{th}$ Jantzen's bilinear form
is being conjugated by $\varphi:\tilde K_m \longrightarrow K_m$ to
the bilinear form
$<\ ,>_{L_\A,f,m}$.\qed

\vskip 5mm\section{The Index of Contact Vector Fields}
 \subsection{The GSV-index $Ind_{V_0,\pm}(X|V_0)$ and the Signature Function
$Sgn_{f,0}(X)$.}

\vskip 5mm Let $f:\R^{2N+1} \rightarrow \R$ be a germ of a real
analytic function with an algebraically isolated singularity at
$0$. Denote also by $f$ its extension to a germ in $\C^{2N+1}$ and
let $V_t$ and $V_t^\C$ be the germs of real or complex analytic
varieties defined by $f=t$. In this section we prove Theorem
\ref{th2}.

We know that both
the GSV-index $Ind_{V_0,\pm}(X|V_0)$ and the signature function
$Sgn_{f,0}(X)$ verify the law of conservation of numbers
(see (\ref{z3}) and similarly for the signature function
$Sgn_{f,0}(X)$
\cite{GM1}). They also coincide in smooth points of the variety.
Hence, the two indices differ by a constant $K_+$ or $K_-$
depending only on the function $f$ (and not on the vector field)
and on the positive or negative sign chosen in the GSV-index.
Given a function $f$ as above, in order to determine these
constants $K_\pm$, it is sufficient to calculate both indices for
one vector field $X_0$ tangent to $V_0$.

\vskip4mm \noindent{\bf Proof of Theorem \ref{th2}:}
In order to prove Theorem~\ref{th2}, we have to study the index of
a family of vector fields tangent to the smoothening $f=t$ of the
singular variety $f=0$. When the ambient space is even
dimensional, this was done (\cite{GM2}) using the Hamiltonian
vector field associated to $f$. Here, we study the odd-dimensional
ambient space $\R^{2N+1}$ and we use the vector fields
$$
 X_t= (f-t)\frac{\partial}{\partial x_1}+ \sum_{i=1}^N[
\frac{\partial f}{\partial x_{2i+1}}\frac{\partial }{\partial
x_{2i}}- \frac{\partial f}{\partial
x_{2i}}\frac{\partial}{\partial x_{2i+1}}],
$$
which we call the contact vector fields. The vector field $X_t$ is
tangent to $V_t$, for any $t$,  since
$D(f-t)X_t=\frac{\partial{f}}{\partial x_1}(f-t)$, where
$\frac{\partial{f}}{\partial x_1}$ is the cofactor. For almost all
linear hyperplanes through $0$ in $\C^{2N+1}$ the projection to
this hyperplane gives a description of $V^\C_0$ as a branched
finite analytic cover \cite{G}. Set $f_j:=\frac{\partial
f}{\partial x_j}$, with $j=1,\ldots,2N+1$. After perhaps a generic
rotation, we may assume that  $0$ is the only point in its
neighborhood that satisfies the equations
$f=f_2=\ldots=f_{2N+1}=0$, or equivalently such that
$f,f_2,\ldots,f_{2N+1}$ is a regular sequence \cite{E}. Hence, the
vector field $X_0$ has an
 algebraically isolated zero at the origin.
 The functions $f_1,\ldots,f_{2N+1}$ form a
regular sequence, since $f$ has isolated singularities. The
hypothesis of the previous part of this paper are satisfied and we
apply Sections 1, 2 and 3 to this situation.
 Choose a small neighborhood $U_\C$  of $0 \in \C^{2N+1}$  and a small
 $\varepsilon>0$,
 as in Section $2.1$.
The derivative of $X_t$  is \begin{equation}\label{3.4} DX_t:=
\left(\begin{matrix} f_1 & f_2 \ \ \  \ \ldots \  \ \ f_{2N+1} \cr
*& \frac{\partial( f_3,-f_2, \ldots,f_{2N+1},-f_{2N})}{\partial
(x_2,\ldots,  x_{2N+1})}
\end{matrix}\right).
\end{equation}
 Denote by
$Y_t:=X_t|_{V_t^\C}$ the restriction of $X_t$ to $V_t^\C$ or to
$V_t$. The singularities of $X_t$ are always contained in
$V_t^\C$, and hence $X_t$ and $Y_t$
 have the same singularities: $\CC_\C \cap V_t^\C$.

\vskip 5mm By definition (see \cite{GSV}), the GSV-index
$Ind_{V_\pm,0}(Y_0)$ is the sum of the indices of $Y_t$ at the
points $p_t\in V_t$, $\pm t>0$ small:
\begin{equation}\label{z3}
Ind_{V_0,\pm}(Y_0,0)= \sum_{
\begin{matrix} p_t \in U \cap V_t, Y_t(p_t)=0 \cr\pm t>0\end{matrix}}
Ind_{V_t}(Y_t,p_t) .\end{equation} Note that $V_t$ is smooth, so
the signatures $Ind_{V_t}(Y_t,p_t)$ can be calculated using the
usual Eisenbud-Levine, Khimshiashvili formula, on the smooth
variety $V_t$. That is, instead of using the Jacobian $J(X_t)$ as
the generator of the socle, one uses the relative Jacobian
$J(Y_t)$. In the localization of the algebra $\B_t$ in $p_t$, we
have
\begin{equation}\label{z7} J(X_t) = f_1
J(Y_t).\end{equation} Hence, the signature of the bilinear form
$< ,\,>^{rel}_t$ (\ref{1.3d_t})
gives the GSV-index:
\begin{equation}\label{GSV}
Ind_{V_\pm,0}(X_0)=\tau_\pm.\end{equation} On the other hand, by
definition (\cite{GM1}), the signature function
$Sgn_{f,0}(X)$
is given by the signature of the form
$< ,\,>^{rel}_t$, for $t=0$. That is,
\begin{equation}
\label{hom1}
Sgn_{f,0}(X_0)
=\tau_0\end{equation}
It now follows from Jantzen's Proposition \ref{Proposition 2.1}
that the constants $K_\pm$ in Theorem \ref{th2} are
\begin{equation}\label{K}
K_+=\sum_{m\geq1}\tau_m,\quad K_-=\sum_{m\geq1}(-1)^m
\tau_m.\end{equation} The Theorem \ref{th2} finally follows from
(\ref{K}) applying Theorem \ref{Proposition 3.1}, which asserts that
$\tau_m=\sigma_{A,f,m}$. \qed
\bigskip

\vskip 3mm \begin{Corollary}\label{Corollary 2} Let $\sigma_\A$
be the signature of the Jacobian algebra $\A$  in (\ref{d3}), and
let $ \sigma_{\A,f,m}$, $m=0,\ldots, \ell+1$, be defined as above.
Then
$$\sigma_A = \frac{\chi_+ - \chi_-}{2} = \sum_{m= odd} \sigma_{\A,f,m}.$$
\end{Corollary}

\noindent{\bf Proof:} By Arnold's formula (\cite{A}),
$2\sigma_\A$ equals ${\chi_+ -\chi_-}$. Now, by the
Poincar\'e-Hopf index theorem, ${\chi_+ -\chi_-}$ equals
$Ind_{V_0,+}(X)-Ind_{V_0,-}(X)$, where $X$ is a real vector field
having an algebraically isolated singularity at the origin tangent
to $V$. The Corollary now follows from Theorem \ref{th2}. \qed

\vskip 5mm \noindent \subsection{Examples} \vskip 3mm \noindent
\begin{Example}\label{Example 4.1} Let $f$ be a quasi-homogeneous
real analytic function with an algebraically isolated singularity,
i.e. $[f]_\A=0 \in \A$. In this case, $Ann_\A(f)=\A$, $M_f=0$ and
$\sigma_{f,1}=\tau_1$ is the only non-zero Jantzen signature of
order higher than 0 and it is equal to $\sigma_A$. Hence $K_\pm =
\chi_\pm=\pm\sigma_\A$.
\end{Example}
\vskip 3mm \noindent
\begin{Example}\label{Example 4.2} Let $f=(x^2+y^3)(x^3+y^2)+z^2$ and
$V=f^{-1}(0)\subset \R^3$. This example is not a quasi-homogeneous
singularity. All calculations have been done using the Computer
algebra system Singular \cite{S}. The local algebra $\A=
\frac{{\mathcal A}_{\R^{3,0}} }{ (f_x,f_y,f_z)}$ has dimension
$11$, $Ann_\A(f)$ is the maximal ideal of dimension 10 and
$([f]_\A)$ is the 1-dimensional socle ideal. We thus have that
$M_f$ has $9$ one-dimensional Jordan blocks and one 2-dimensional
Jordan block. The Hessian $[Hess(f)]_\A$ generating the socle
equals $-220[f]_\A$ in $\A$. The filtration (\ref{a}) is given by
$$
(f)\subset Ann_\A(f)\subset\A.$$
 The signature $\sigma_1$ is the signature of
$< ,>_{L_\A}$ on the $9$ dimensional space isomorphic to
$\frac{Ann_\A(f)}{(f)}$. This signature is equal to $1$. The
signature $\sigma_2$ is given by the the sign of
$L_\A(\frac{f\cdot f}{f})=L_\A(f)<0,$ so $\sigma_2$=-1. This gives
by Theorem~2 that $K_+=0$, $K_-=-2$
\end{Example}


\begin{thebibliography}{99}
{\small}
\bibitem{A} V. Arnold, Index of a singular point of a vector field, the
Petrovskii-Oleinik inequality and mixed Hodge structures,
Functional Analysis and its Applications, Vol 12, No 1, (1978),
1--14.
\bibitem{BG} Ch. Bonatti, X. G\'omez-Mont, The index of a holomorphic vector
field on a singular variety, I, Ast\'erisque, 222 (1994), 9--35.
\bibitem{EL} D. Eisenbud, H. Levine, An algebraic formula for the degree of a
$C^\infty$ map germ, Ann. Math., 106 (1977), 19--38.
\bibitem{E} D. Eisenbud, Commutative Algebra with a view towards Algebraic
Geometry, Springer 1995.
\bibitem{GGM} L. Giraldo, X. G\'omez-Mont, P. Marde\v si\'c
Computation of topological numbers via linear algebra:
Hypersurfaces, vector fields and vector fields on hypersurfaces,
Contemporary Math. Vol 240 (1999), 175--182.
\bibitem{GM1} X. G\'omez-Mont, P. Marde\v si\'c, The index of a vector field
tangent to a hypersurface and the signature of the relative
Jacobian determinant, Annals Inst. Fourier 47, 5 (1997), 1523 --
1539.
\bibitem{GM2} X. G\'omez-Mont, P. Marde\v si\'c, The index of a vector field
tangent to an odd dimensional hypersurface and the signature of
the relative Hessian, Functional Analysis and its Applications,
Vol 33, No 1, (1999).
\bibitem{GSV} X. G\'omez-Mont, J. Seade, A. Verjovsky, The index of a
holomorphic flow with an isolated singularity, Math. Ann., 291
(1991), 737--751.
\bibitem{S} G. -M. Greuel, G. Pfister. H. Sch\"onemann, {\sc Singular}
3.0. A Computer Algebra System for Polynomial Computations. Center
for Computer Algebra, University of Kaiserslautern (2005). {\tt
http://www.singular.uni-kl.de}.
\bibitem{G} R. C. Gunning, Introduction to Holomorphic Functions of Several Complex
Variables, Vol. II: Local theory, Wadsworth \& Brooks/Cole,
(1990).
\bibitem{J} J.C. Jantzen,  Moduln mit einem H\"ochsten Gewicht, Lecture
Notes in Mathematics 750, Springer, Berlin-New York- Tokio,
(1979).
\bibitem{K} G. N. Khimshiashvili, On the local degree of a smooth map. Trudi Tbilisi math. Inst, (1980), 105-124.
\bibitem{MvS} J. Montaldi, D. van Straten, One-forms on singular curves and the topology
of real curves singularities, Topology 29, no.4 (1990), 501--510
\bibitem{Vo} D. Vogan, Unitarizability of certain series of representations,
Ann. of Math.(2) 120 (1984) no.1. 141-187.

\end{thebibliography}
\end{document}